\documentclass[a4paper,USenglish,cleveref,thm-restate,preprint]{lipics-v2021}
\title{A Vanishing Lower Ratio for Primes of the Form $n^2+s^2$}
\titlerunning{A Vanishing Lower Ratio for $n^2+s^2$}
\author{Sheikh Azizul Hakim}{Carnegie Mellon University, Pittsburgh, Pennsylvania, USA}{}{}{}
\authorrunning{S. A. Hakim}
\Copyright{Sheikh Azizul Hakim}
\relatedversion{}
\EventEditors{}
\EventNoEds{0}
\EventLongTitle{}
\EventShortTitle{}
\EventAcronym{}
\EventYear{2026}
\EventDate{}
\EventLocation{}
\EventLogo{}
\SeriesVolume{0}
\ArticleNo{0}
\begin{document}
\nolinenumbers
\maketitle
\begin{abstract}
For an integer $n\ge 2$, let $a(n)$ count the integers $s$ with
$1\le s<n$ for which $n^2+s^2$ is prime. Let $\pi(n)$ count the primes at
most $n$. The lower bound $a(n)\ge\pi(n)/5$ was conjectured in OEIS A069004.
A proof is given that the ratio $a(n)/\pi(n)$ comes arbitrarily close to
zero for infinitely many values of $n$. The construction chooses $n$
divisible by many primes congruent to $3$ modulo $4$. The associated
large-sieve estimate yields the claimed result through an Euler product,
Karamata's Tauberian theorem, and the prime number theorem.
\end{abstract}

\section{The conjecture and the construction}
For each integer $n\ge2$, let $a(n)$ be the number of integers $s$
between $1$ and $n-1$ for which $n^2+s^2$ is prime. Let $\pi(n)$ be the
number of primes at most $n$. Formally, using $\mathbb Z$ for the integers
and $\#$ for the number of elements in a finite set, these definitions are
\[
 a(n)=\#\{s\in\mathbb Z\mid 1\le s<n,\ n^2+s^2\text{ is prime}\},
 \qquad \pi(n)=\#\{p\le n\mid p\text{ is prime}\}.
\]
For example, $a(5)=2$, because $5^2+2^2=29$ and $5^2+4^2=41$ are prime.
The lower-bound conjecture attributed to T.~D. Noe in OEIS A069004
\cite{oeis} asserted that $5a(n)\ge\pi(n)$ for every $n>1$.

The elementary divisibility step does not require a congruence condition.
For \emph{any} prime $p$ dividing $n$, one has
\[
 n^2+s^2\equiv s^2\pmod p,
 \qquad
 p\mid(n^2+s^2)\iff p\mid s.
\]
Indeed, a prime divides $s^2$ if and only if it divides $s$.
The reason to select primes $p\equiv3\pmod4$ is different. When such a
prime does not divide $n$, the congruence $s^2\equiv-n^2\pmod p$
has no solution, since $-1$ is not a square modulo $p$. Making $p$ divide
$n$ therefore \emph{creates} one forbidden residue class, namely $s\equiv0$.
By contrast, a prime $p\equiv1\pmod4$ coprime to $n$ already excludes
two classes, and making it divide $n$ would reduce that number to one.
Hence primes congruent to $3$ modulo $4$ improve this particular sieve,
whereas primes congruent to $1$ modulo $4$ do not. The large sieve counts
integers avoiding all the forbidden residue classes simultaneously.

Let $\mathcal P$ be any finite set of primes congruent to $3$
modulo $4$. Define $M$ as the product of its primes and $B$ as the
product of the corresponding factors $p/(p-1)$. For every positive
integer $r$, define $n_r$ by
\[
 M=\prod_{p\in\mathcal P}p,\qquad
 B=\prod_{p\in\mathcal P}\frac{p}{p-1},\qquad
 n_r=2^rM.
\]
An empty product equals $1$. In limits as $r$ tends to infinity,
$\mathcal P$, $M$, and $B$ remain fixed. The notation $\limsup$ and
$\liminf$ refers to the upper and lower limits of a sequence.

\begin{theorem}\label{thm:main}
For every finite set $\mathcal P$ as above,
\[
 \limsup_{r\to\infty}\frac{a(n_r)}{\pi(n_r)}
 \le\frac{2\pi}{B}.
\]
Consequently,
\[
 \boxed{\displaystyle\liminf_{n\to\infty}\frac{a(n)}{\pi(n)}=0.}
\]
In particular, infinitely many integers $n>1$ satisfy
$5a(n)<\pi(n)$.
\end{theorem}

\section{Which residue classes are excluded?}
For a prime $p$ and a positive integer $n$, let $\Omega_p(n)$
be the set of residue classes modulo $p$ that satisfy
$s^2\equiv -n^2\pmod p$. In symbols,
\[
 \Omega_p(n)=\{s\bmod p\mid s^2\equiv-n^2\pmod p\}.
\]
Fix the set $\mathcal P$ chosen above. When $n=n_r$, the number of
classes in $\Omega_p(n_r)$ does not depend on $r$. Denote that number
by $\omega_{\mathcal P}(p)$, where $|S|$ means the number of elements of
any finite set $S$. Thus
\[
 \omega_{\mathcal P}(p)=|\Omega_p(n_r)|.
\]
The residue-class count is
\begin{equation}\label{eq:omega}
 \omega_{\mathcal P}(p)=
 \begin{cases}
 1,&p=2,\\
 2,&p\equiv1\pmod4,\\
 1,&p\in\mathcal P,\\
 0,&p\equiv3\pmod4\text{ and }p\notin\mathcal P.
 \end{cases}
\end{equation}
Indeed, modulo $2$ the congruence has one solution. Every prime
congruent to $1$ modulo $4$ is coprime to $n=2^rM$, and $-1$ has exactly
two square roots modulo such a prime. For $p\equiv3\pmod4$, $-1$
is a non-square. Hence if $p\nmid n$ there are no solutions, whereas
if $p\mid n$ the unique solution is $s\equiv0\pmod p$. The latter
conclusion uses only the elementary divisibility equivalence above, not
the congruence condition on $p$.
In all cases $0\le\omega_{\mathcal P}(p)<p$.

Define an arithmetic function $h_{\mathcal P}$ by setting
$h_{\mathcal P}(1)=1$, setting $h_{\mathcal P}(p^j)=0$ for each prime
$p$ and integer $j\ge2$, and multiplying these values across
coprime factors. Its value at a prime $p$ is
\begin{equation}\label{eq:h}
 h_{\mathcal P}(p)=\frac{\omega_{\mathcal P}(p)}{p-\omega_{\mathcal P}(p)}=
 \begin{cases}
 1,&p=2,\\
 \dfrac{2}{p-2},&p\equiv1\pmod4,\\
 \dfrac1{p-1},&p\in\mathcal P,\\
 0,&p\equiv3\pmod4,\ p\notin\mathcal P.
 \end{cases}
\end{equation}
For each real number $x\ge1$, define the cumulative sum $H_{\mathcal P}(x)$ by
\[
 H_{\mathcal P}(x)=\sum_{1\le d\le x}h_{\mathcal P}(d),
\]
where $d$ runs over positive integers.

\begin{lemma}[Sieve inequality]\label{lem:sieve}
Let $n=n_r$ and let $Q=\lfloor\sqrt n\rfloor$, the largest integer
not exceeding $\sqrt n$. Then
\[
 a(n)\le\frac{n-1+Q^2}{H_{\mathcal P}(Q)}
 \le\frac{2n}{H_{\mathcal P}(Q)}.
\]
\end{lemma}
\begin{proof}
If $n^2+s^2$ is prime, then it exceeds $Q$, so it is not divisible by
any prime $p\le Q$. Therefore $s\bmod p\notin\Omega_p(n)$ for every such
prime. The arithmetic large sieve in its standard sifted-set form
\cite[Theorem~16.4]{kedlaya} applies to an interval containing $N$
consecutive integers and a collection of excluded residue classes
modulo each prime $p\le Q$. Let $\mu$ denote the M\"obius function,
which is zero on integers divisible by a square and is $(-1)^k$ on
a product of $k$ distinct primes. In the present setting the large
sieve bounds the number of integers avoiding the excluded classes by
\[
 \frac{N+Q^2}{\displaystyle\sum_{1\le d\le Q}\mu^2(d)
 \prod_{p\mid d}\frac{\omega_{\mathcal P}(p)}{p-\omega_{\mathcal P}(p)}}.
\]
Apply this with the interval $1\le s<n$, of length $N=n-1$, and use
\eqref{eq:h} to identify the denominator with $H_{\mathcal P}(Q)$.
Finally, $Q^2\le n$.
\end{proof}

\section{Why the sieve denominator grows}\label{sec:denominator}
The central point is that the normal large-sieve denominator grows
proportionally to $\log Q$, while each selected prime $p\equiv3\pmod4$
multiplies its leading constant by $p/(p-1)$.

Define the numerical constant $\kappa_0$ by
\begin{equation}\label{eq:kappa}
 \kappa_0=\frac\pi4
 \prod_{p\equiv1\ (4)}\left(1+\frac1{p(p-2)}\right)
 \prod_{p\equiv3\ (4)}\left(1-\frac1{p^2}\right).
\end{equation}
Here each infinite product runs over primes in the indicated residue
class. Their convergence and the lower bound $\kappa_0\ge2/\pi$
are established in the proof below. We write $f(x)\sim g(x)$ if
$f(x)/g(x)$ tends to $1$ in the indicated limit.

\begin{lemma}[Denominator asymptotic]\label{lem:denominator}
For any fixed finite set $\mathcal P$,
\begin{equation}\label{eq:den}
 H_{\mathcal P}(x)\sim\kappa_0 B\log x
 \qquad(x\to\infty).
\end{equation}
Moreover, $\kappa_0\ge2/\pi$.
\end{lemma}
\begin{proof}
Let $\chi_4$ be the real Dirichlet character modulo $4$. It
vanishes on even integers, is $1$ on integers congruent to $1$
modulo $4$, and is $-1$ on integers congruent to $3$ modulo $4$.
Write $\zeta(z)$ for the Riemann zeta function and $L(z,\chi_4)$
for the Dirichlet $L$-function of this character, where $z$ is a
complex variable.

First take $\mathcal P=\varnothing$, the empty set. For a positive
real number $t$, define the Dirichlet series $D_\varnothing(t)$ by
\[
 D_\varnothing(t)=\sum_{d\ge1}\frac{h_\varnothing(d)}{d^t}
 =(1+2^{-t})\prod_{p\equiv1\ (4)}
 \left(1+\frac{2p^{-t}}{p-2}\right).
\]
The equality follows from the multiplicative definition of
$h_\varnothing$, and both expressions converge for $t>0$.
For $t\ge0$, define an auxiliary Euler product $K(t)$ by
\begin{align*}
 K(t)={}&(1+2^{-t})(1-2^{-1-t})\\
 &\cdot\prod_{p\equiv1\ (4)}
 \left(1+\frac{2p^{-t}}{p-2}\right)(1-p^{-1-t})^2\\
 &\cdot\prod_{p\equiv3\ (4)}(1-p^{-2-2t}).
\end{align*}
Comparing Euler factors now gives
\begin{equation}\label{eq:factor}
 D_\varnothing(t)=\zeta(1+t)L(1+t,\chi_4)K(t)
 \qquad(t>0).
\end{equation}
For $0\le t\le1$, each factor in the two infinite products defining
$K(t)$ differs from $1$ by at most a constant times $p^{-2}$.
It follows that the products converge uniformly on this interval and
$K(t)\to K(0)>0$ as $t\downarrow0$.
Evaluation at zero gives
\[
 K(0)=\prod_{p\equiv1\ (4)}\frac{(p-1)^2}{p(p-2)}
 \prod_{p\equiv3\ (4)}\left(1-\frac1{p^2}\right).
\]
Since $\zeta(1+t)\sim 1/t$ and
$L(1,\chi_4)=1-1/3+1/5-\cdots=\pi/4$,
\[
 D_\varnothing(t)\sim\frac{\kappa_0}{t},
 \qquad\kappa_0=\frac\pi4K(0).
\]
For general fixed $\mathcal P$, exactly the Euler factors with
$p\in\mathcal P$ are inserted. Hence
\[
 \sum_{d\ge1}\frac{h_{\mathcal P}(d)}{d^t}
 =D_\varnothing(t)
 \prod_{p\in\mathcal P}\left(1+\frac{p^{-t}}{p-1}\right)
 \sim\frac{\kappa_0B}{t}.
\]
To convert this into a partial-sum asymptotic, let $v\ge0$ be
real and define the nondecreasing function
\[
 U(v)=\sum_{\substack{d\ge1\\ \log d\le v}}h_{\mathcal P}(d).
\]
Its Laplace--Stieltjes transform is the Dirichlet series displayed
above. Karamata's Tauberian theorem in the case of index $1$
\cite[Theorem~1.7.1]{bgt} yields $U(v)\sim\kappa_0Bv$ as $v\to\infty$.
Taking $v=\log x$ proves \eqref{eq:den}.

Finally, every factor in the first product of \eqref{eq:kappa}
is at least $1$, while the second product is at least the product over
\emph{all} odd primes. Euler's evaluation $\zeta(2)=\pi^2/6$ gives
\[
 \prod_{p\equiv3\ (4)}\left(1-\frac1{p^2}\right)
 \ge\prod_{p>2}\left(1-\frac1{p^2}\right)
 =\frac{8}{\pi^2}.
\]
Thus $\kappa_0\ge(\pi/4)(8/\pi^2)=2/\pi$.
\end{proof}

\section{Proof of the theorem}
\begin{proof}[Proof of Theorem~\ref{thm:main}]
Fix $\mathcal P$, put $n=n_r$, and let $Q=\lfloor\sqrt n\rfloor$.
As $r\to\infty$, we have $n\to\infty$ and
$\log Q\sim\frac12\log n$. In the next display, $o(1)$ denotes
a quantity tending to zero as $r\to\infty$.
Lemmas~\ref{lem:sieve} and \ref{lem:denominator} imply
\[
 a(n)\le\frac{2n}{(\kappa_0 B+o(1))\log Q}
 =\left(\frac4{\kappa_0 B}+o(1)\right)
 \frac{n}{\log n}.
\]
By the prime number theorem, $\pi(n)\sim n/\log n$.
Consequently
\[
 \limsup_{r\to\infty}\frac{a(n_r)}{\pi(n_r)}
 \le\frac4{\kappa_0 B}\le\frac{2\pi}{B}.
\]
It remains to show that $B$ can be made arbitrarily large.
Let $u$ be a real number greater than $1$. The Euler products of
$\zeta(u)$ and $L(u,\chi_4)$ give the following identity up to a
bounded error. Here $O(1)$ denotes a quantity that stays bounded as
$u$ decreases to $1$.
\[
 \sum_{p\equiv3\ (4)}\frac1{p^u}
 =\frac12\bigl(\log\zeta(u)-\log L(u,\chi_4)\bigr)+O(1)
 \longrightarrow\infty\qquad(u\downarrow1),
\]
because $L(1,\chi_4)=\pi/4>0$. Hence
$\sum_{p\equiv3\ (4)}1/p$ diverges, and so does
$\sum_{p\equiv3\ (4)}\log(p/(p-1))$.
Thus finite sets $\mathcal P$ can be chosen with $B$ arbitrarily large.
For any $\varepsilon>0$, choose $\mathcal P$ with
$2\pi/B<\varepsilon$. The preceding limsup bound then shows that
$a(n_r)/\pi(n_r)<\varepsilon$ for all sufficiently large $r$.
Since these $n_r$ tend to infinity, the asserted limit inferior is zero.
Taking $B>10\pi$ makes the limsup strictly smaller than $1/5$,
so all sufficiently large $r$ yield the desired counterexamples.
\end{proof}

\section{Discussion and verification}

The proof presented here was inspired by a Lean 4 formalization
developed through the \emph{Left to Solve} initiative
\cite{lefttosolve}. The original formal proof is publicly
available on GitHub \cite{leanproof}, together with its
submission to the Open Math Agents project \cite{pr}.

The central construction, in which $n=2^rM$ and $M$ is
a product of selected primes congruent to $3$ modulo $4$,
is motivated by the Lean proof. However, the present
argument takes a different analytic approach. Rather
than deriving explicit finite bounds through Selberg's
sieve, it applies the large sieve together with an
Euler product and a Tauberian theorem.

This yields a stronger asymptotic conclusion. The ratio
$a(n)/\pi(n)$ has limit inferior zero, so infinitely
many counterexamples exist. For the original conjecture,
a choice of $B>40$ suffices, whereas the Lean proof
uses the larger threshold $200000$ and specifies an
explicit exponent. The present argument establishes
existence for all sufficiently large exponents without
giving a numerical threshold.

The Lean submission has passed an automated proof check
\cite{pr}. The mathematical proof given here is inspired
by that formalization but is not itself formalized in
Lean. Its additional asymptotic conclusion therefore
requires independent mathematical verification.

\paragraph*{AI-use disclosure.}
AI-based systems (Gemini (as in Antigravity), ChatGPT, and Claude) contributed to discovering and formalizing a counterexample
strategy, developing the alternative large-sieve exposition presented here,
and drafting the manuscript. The author is responsible for evaluating the
mathematical argument, the attribution, and the final submitted version.


\begin{thebibliography}{9}
\bibitem{oeis}
OEIS Foundation Inc.
\emph{OEIS A069004, Number of times $n^2+s^2$ is prime for $0<s<n$}.
On-Line Encyclopedia of Integer Sequences.
\url{https://oeis.org/A069004}.
Conjecture attributed to T. D. Noe.



\bibitem{kedlaya}
Kiran S. Kedlaya.
\emph{An Introduction to Analytic Number Theory}, Chapter 16, Theorem 16.4.
Online lecture notes. \url{https://kskedlaya.org/ant/chap-largesieve2.html}.

\bibitem{bgt}
N. H. Bingham, C. M. Goldie, and J. L. Teugels.
\emph{Regular Variation}.
Cambridge University Press, 1987. Theorem 1.7.1.

\bibitem{pr}
skhakim.
Disproof of OEIS A069004 conjecture 2 without upper bound.
Pull request 2077, \emph{open-math-agents}, 2026.
\url{https://github.com/GuillaumeBAUD-coder/open-math-agents/pull/2077}.


\bibitem{lefttosolve}
Left to Solve.
\emph{Open Mathematical Problems and AI-Assisted
Formal Proofs}.
\url{https://lefttosolve.com/}.

\bibitem{leanproof}
Open Math Agents.
\emph{Lean 4 Formalization of a Disproof of
OEIS A069004, Conjecture 2}.
\url{https://github.com/skhakim/open-math-agents/blob/oma/OeisA69004.conjecture2.variants.without_upper_bound/proofs/OeisA69004.conjecture2.variants.without_upper_bound.lean}.

\end{thebibliography}
\end{document}